\documentclass[11 pt]{amsart}
\usepackage{plain, times, txfonts}
\usepackage{graphicx, mathrsfs}
\usepackage{amssymb,amsmath,amsthm}
\usepackage{mathrsfs,dsfont}

\theoremstyle{plain}
\newtheorem{theorem}{Theorem}[section]
\newtheorem*{nth}{Theorem}
\newtheorem{proposition}{Proposition}[section]

\numberwithin{equation}{section} 

\newcommand\sms{Schr\"odinger maps}
\newcommand\sm{Schr\"odinger map}

\newcommand\intl{\int\limits}
\newcommand\<{\langle}
\renewcommand\>{\rangle}

\def\CD{{\mathcal D}}

\def\CL{{\mathcal L}}

\def\CC{\Bbb{C}}
\def\HH{{\mathbb H}}
\def\RR{{\mathbb R}}
\def\R{{\mathbb R}}
\def\SS{{\mathbb S}}
\def\UU{{\Bbb U}}

\def \te{{\tilde e}}
\newcommand{\lf}{\left}
\newcommand{\rt}{\right}
\newcommand{\p}{\partial}
\newcommand{\half}{\frac{1}{2}}
\newcommand{\bfi}[1]{\textbf{\textit{#1}}}

\newcommand{\sdef}{\stackrel{\rm def}{=}}

\newcommand{\he}{\hat{e}}
\newcommand{\ta}{\tilde{a}}

\newcommand{\hA}{\hat{A}}

\newcommand{\hq}{\hat{q}}
\newcommand{\tq}{\tilde{q}}
\newcommand{\tu}{\tilde{u}}

\newcommand{\ha}{\hat{a}}
\newcommand{\hf}{\hat{f}}

\newcommand\h{\mathbb{H}}
\newcommand\op{\operatorname}

\theoremstyle{remark}

\begin{document}

\title[Schr\"odinger Maps and Frame Systems ]{Schr\"odinger Maps and their associated Frame Systems}
\author[Nahmod]{Andrea Nahmod$^1$}
\address{$^1$Department of Mathematics \\
University of Massachusetts\\ 710 N. Pleasant Street, Amherst MA 01003}
\email{nahmod@math.umass.edu}
\thanks{$^1$ The first author is funded in part by NSF DMS 0503542.}
\author[Shatah]{Jalal Shatah$^2$}
\thanks{$^2$ The second author is funded in part by NSF DMS
0203485.}
\address{$^2$Courant Institute of Mathematical Sciences\\
251 Mercer Street\\ New York, NY 10012} \email{shatah@cims.nyu.edu}
\author[Vega]{Luis Vega$^3$}
\address{$^3$ Universidad del Pa\'is Vasco, Departamento de Matem\'aticas, Apto. 644, 48080 Bilbao, Spain}
\email{mtpvegol@lg.ehu.es}
\thanks{$^3$ The third author is funded in part by MTM 2004-03029 of MEC (Spain) and FEDER.}
\author[Zeng]{Chongchun Zeng$^4$}
\address{$^4$School of Mathematics\\
Georgia Institute of Technology\\ Atlanta, GA 30332}
\email{zengch@math.gatech.edu}
\thanks{$^4$ The fourth author is funded in part by NSF
DMS 0627842 and the Sloan Fellowship.}
\thanks{ Part of this work was done while the first and third authors were members at the Institute for Advanced Study and the fourth was visiting the Courant Institute}
\date{}
\begin{abstract}
In this paper we establish the equivalence of solutions between \sms\ into $\SS^2$ or $ \HH^2$ and their associated gauge invariant Schr\"odinger equations.  We also establish the existence of global weak solutions into $\HH^2$ in two space dimensions. We extend these ideas for maps into compact hermitian symmetric manifolds with trivial first cohomology.
\end{abstract}
\maketitle

\section{Introduction}
Schr\"odinger maps are maps from space-time into  a K\"ahler manifold with metric
$h$ and complex structure $J$ satisfying:  $u : \RR^d \times \RR \to (M, h, J)$
\begin{equation}
\partial_t u =J \sum_\ell D_\ell \partial^\ell u, \tag{SM}
\end{equation}
where $D$ denotes the covariant derivative on $u^{-1}TM$.
These  maps are a generalization of  the Heisenberg model describing the magnetization spin $m \in \SS^2 \subset \RR^3$ in a ferromagnetic material
\begin{equation*}
\partial_t m = m \times \Delta m.
\end{equation*}
For $m \in \SS^2$  the operator $J= m\times$
acting on $T_m \SS^2$ is equivalent to complex multiplication by
$i$ on $\CC$.  Thus
\[
\p_tm =  m \times \Delta m = m\times ( \Delta m + | \nabla m |^2 m ) =  J\sum_i D_i \p^i m
\]
where as before $D_i= \p_i\phantom{v} +\< \p_im, \phantom{v}\> m$ denotes the covariant derivative on $m^{-1}T\SS^2$; and the Heisenberg model  can be written as
\[ \partial_t m = J \sum_iD_i \p^i m .
\]

In one-space dimension the Heisenberg model can be
transformed into the focusing NLS
\begin{equation}\tag{NLS}
i \p_t q - \p^2_x q - \half   |q|^2   q= 0
\end{equation}
via the Hasimoto transformation.  This transformation was later
generalized by N. H. Chang, J. Shatah, and K. Uhlenbeck \cite{CSU} to study the regularity
of Schr\"{o}dinger maps. The idea  in \cite{CSU} was
to disregard  the customary coordinates representation of the (SM) system  and to
introduce instead  a gauge invariant nonlinear Schr\"odinger equations (GNLS) derived by using a pull-back frame on $u^{-1} TM$.  The GNLS is given schematically  by
$$
\CD_t  q = i   \sum_k \CD^2_k   q  + i   F q.
$$
Using the Coulomb gauge, this sytem can be written as
$$
i\p_tq = \Delta q + \Delta^{-1}[\p(\mathcal{O}(|q|^2)]\p q + \mathcal{O}(|q|^3).
$$
One of the consequences of such a representation  was to reveal the semilinear nature of the Schr\"{o}dinger
maps systems which led to  the first regularity proof in 1 and 2-dimensions for finite energy equivariant data \cite{CSU}.  Here we would like to note that the  $1$-dimensional Cauchy problem for (SM) is subcritical with respect to the energy space $\p u \in L^2$ and as such should be solvable for data $\p u\in L^2$.  However the only proof of global well-posedness in this case was given for data $\p  u \in H^1$ and uses the GNLS system \cite{CSU}.  The desired goal would be to solve the Cauchy problem  and to show equivalence  when the derivative of the data  behaves like $\delta(x)$; i.e. data scaling as $\p u \in {\dot{H}}^{-1/2}$.

Another consequence of introducing the GNLS was to show that  for constant curvature  $M$ the GNLS system doesn't depend explicitly on $u$, and therefore can be solved  without any reference to the SM system.

Using this last observation a natural question to ask in the constant curvature case  is:  When do solutions of the GNLS represent solutions of SM? \,
For smooth solutions this question was answered in one dimension by Terng and Uhlenbeck \cite{TEU} and in two dimensions, for a special case, by N.H. Chang and O. Pashev \cite{CHPSV}.

In this paper we are interested in studying the correspondence between solutions $u$ of  the Schr\"odinger map system  and solutions $q$ of its associated  gauge invariant nonlinear Schr\"odinger equations for low regularity data.  In particular we show the equivalence of the two systems for solutions where the problems are expected to be well posed, i.e., $\p u\in H^{\frac d2-1}$ plus Strichartz  estimates for $d = 2$. One should remark that the interesting cases for the equivalence of the SM system and  GNLS system correspond to $d=1, 2$ or    $3$ since  in $d\ge4$,  $\p u\in H^{\frac d2-1}$ and  equation (SM) holds a.e.; thus there is little difference between smooth and $\p u \in H^{\frac d2-1}$  solutions.

The outline of the paper is as follows: In section \ref{frame} we present the frame system.  In section \ref{sphere}  we study the equivalence problem when the target is the sphere. For $d=1$ we show the equivalence under the condition $q\in L^2_{tx}$ and $|q|^2\in L^2(H^{-1})$. For higher dimension we show the equivalence of GNLS and SM  for data in critical spaces, i.e., invariant under the scaling $u(t,x)\to u(\lambda^2 t, \lambda x)$.   In section \ref{hyperbolic} we study the problem when the target is $\HH^2$, the hyperbolic space.  Here we  show equivalence of smooth solutions and for two space dimensions we show global existence of finite energy solutions. Finally in section
\ref{Epilogue} we describe the extension of these results for maps into compact hermitian symmetric manifolds with trivial first cohomology.
\smallskip

Throughout  this paper we sum over repeated indices unless we explicitly state the contrary, and we follow the convention that Greek subscripts vary from $0$ to $d$ while roman subscript vary from $1$ to $d$ or $n$ depending on the context.

\section{Frame System}\label{frame}
The use of frames on the pullback bundle was introduced in \cite{CSU}, and was later used successfully  to study the Cauchy problem for wave maps \cite{SS02, NSU03b}. In \cite{NSU03a} similar ideas as in \cite{CSU} were also used, starting with the pull-back of the conformal frame of $\Bbb S^2$ -which amounts to the stereographic projection- followed by the Coulomb gauge transformation.

\subsubsection*{\bfi{Frames on the pullback bundle}} Let $\phi: \RR^d \to (M, h, J)$ be a map into a  $2n$-dimensional K\"ahler manifold and
let $D$ denote the covariant derivative on $\phi^{-1}   TM$. Since
$M$ is K\"ahler then $D_\ell J(\phi(x))=0$ for $\ell=1, \dots, d$.  With
a slight abuse of language we will refer to sections on $\phi^{-1}  TM$
as vectors.  Let $\{e_a\}^{2n}_{a=1}$ denote an orthonormal frame on
$\phi^{-1}  TM$ such that $e_{a+n}=J  e_a$ for $a=1, \dots, n$.  Such
a frame always exists since $\RR^d$ is contractible  and $M$ is K\"ahler.
\begin{proposition}\label{global-f}
Fix the origin $0\in \RR^d$ and introduce polar coordinates $(r,\omega)$ on $\RR^d$.  Given a  smooth $\phi: \RR^d \to (M, h, J)$, let  $\{e^*_1, \dots e^*_{2n}\} $ denote an orthonormal set of vectors on $\phi^{-1}   TM$ at $x=0$ such that $e^*_{a+n}=J  e^*_a $  and let   $\{e_1, \dots e_{2n}\}$ be the solution to the ODEs
$$
D_re_a =0, \quad
e_a(0,\omega) = e^*_a.
$$
Then $\{e_1, \dots e_{2n}\}$  is an orthonormal frame for   $\phi^{-1}   TM$ with $e_{a+n}=J  e_a$ for $a=1, \dots, n$.
\end{proposition}
\begin{proof}Solve the linear ODEs and use the fact that $D_rJe_a=JD_re_a$  since $M$ is K\"ahler. \end{proof}
Write the frame as $\{e_1, \dots e_{2n}\} = \{e_1, \dots e_n,
Je_1, \dots, Je_n\}\sdef \{e,Je\}$.  For any vector  $v\in \phi^{-1}  TM$ with coordinates $v=\sum^{2n}_{\ell=1} v_\ell    e_\ell$, we introduce complex coordinates $w = (w_1, \cdots , w_n)\in\CC^n$, where $w_\ell =v_\ell +i  v_{\ell+n}$, on $\phi^{-1}  TM$ and write $v=w\cdot e$ where
$$
v=\sum^{2n}_{\ell=1} v_\ell    e_\ell = \sum^n_{\ell=1}  \big(v_\ell + v_{\ell+n}J\big)
e_\ell =   \sum^n_{\ell=1}  \big(v_\ell +i  v_{\ell+n}\big)
e_\ell= \sum^n_{\ell=1} w_\ell   e_\ell = w\cdot e.
$$
In these complex coordinates, $J\to i$ on  $\phi^{-1}TM$.

The covariant derivative on $M$ introduces a connection $\{A_\ell \}$ on $\phi^{-1}  TM$ given by $D_\ell  e_a= A^b_{a \ell}\,  e_b$ for $a=1, \dots, n, \quad \ell=1, \dots, d$.  We simply write
\begin{equation}\label{gauge}
D_\ell  e =  A_\ell\cdot e.
\end{equation}
where the   $n\times n$ matrices   $A_\ell=\big(A^b_{a\ell}\big)\in\mathfrak{su}(n)$.  For  any vector $v= w\cdot e\in \phi^{-1}  TM$, with coordinates $w\in \CC^n$ we have
\begin{equation*}
D_\ell  v = D_\ell (w \cdot e) =  \big(\p_\ell  w + A_\ell  w\big) \cdot e \sdef(\CD_\ell  w) \cdot e
\end{equation*}
where  $\CD$ denotes the covariant derivative on  $\phi^{-1}TM$ expressed in terms of the frame $\{e, Je\}$.

If one chooses another frame $\{\he,J\he\}$ related to $\{e, Je\}$ by a transformation $g\in \SS\UU(n)$, i.e., $\he=g\cdot e$ then
$$\begin{aligned}
&D_\ell  \he = \hA_\ell\cdot \he\\
&\hA_\ell =g^{-1}A_\ell g + g^{-1}\p_\ell g.
\end{aligned}
$$
Thus fixing a frame is equivalent to fixing the connection $A $; {\it i.e.}  fixing a gauge. The matrices $\{A_\ell\}_{\ell =0}^d$, which are given by \eqref{gauge},  have to verify the curvature
equation. That is,  if we let $\p_k\phi = q_k\cdot e$ and denote by
\begin{equation*}\label{eq:rcurv}
\lf[D_\ell, D_k\rt] e_a = R\big(\p_k    \phi, \p_\ell   \phi\big)e_a = R\big(q_k\cdot e, q_\ell \cdot e\big)e_a \sdef F(q_ \ell, q_k)\cdot e_a = F_{ \ell k}\cdot e_a
\end{equation*}
then  we have
\begin{equation}\label{eq:curv}
\lf[\CD_\ell, \CD_k\rt] = \p_\ell   A_k -\p_k    A_\ell
+ \lf[A_\ell, A_k\rt]
= F_{\ell  k}.
\end{equation}
Here it is worth mentioning that the frame constructed in proposition \ref{global-f} corresponds to choosing a connection such that  $x^kA_k(x)=0$.  This gauge is referred to as the exponential (or Cr\"omstrom)  gauge \cite{Uh83}.  For this gauge the connection $A$ can be easily recovered from $F$ by the formula
\begin{equation}\label{exp}
A_k(x) = \int_0^1 x^\ell F_{\ell k}(sx)\, s ds.
\end{equation}

Throughout this paper we are interested in a special frame which corresponds to the  Coulomb gauge, i.e., a frame for which $\sum_1^d \p_\ell \hA_\ell =0$.   Local smooth Coulomb frames can always be constructed as was demonstrated by K. Uhlenbeck in \cite{Uh83}.  This is done by solving the elliptic equation for $g$
$$
0=\p_\ell\hA_\ell = \p_\ell\big(g^{-1}A_\ell g + g^{-1}\p_\ell g\big),
$$
locally on balls in $\RR^d$.  For $d>1$ gluing these local solutions does not necessarily  yield a global Coulomb frame.  Of course if $n=1$ then $g=\exp(i\theta)$ and the above equation is linear and can be solved globally.  For the general problem Dell Antonio and Zwanziger  \cite{DELLZ}  showed that the  existence of a global  ${\dot{ H}}^1$ Coulomb frame.
 \begin{proposition}\label{pro:co}
Given a smooth map $\phi : \RR^d \to (M, h, J) $\,  there exists a frame $\{\he, J\he\}$  such that
\begin{align*}
&D _\ell   \he_a = \hA _\ell \cdot e\\
&\sum^d_{\ell=1} \p _\ell   \hA _\ell =0.
\end{align*}
\end{proposition}
\noindent {\it Sketch of the proof.}  Fix a  frame $\{e,Je\}$  of $\phi^{-1}TM$ and let $A_\ell$ be given by $D_\ell e = A_\ell\cdot e$.  For any $g\in \SS\UU(n)$ let  $\hA_\ell =g^{-1}A_\ell g + g^{-1}\p_\ell g$ and consider the variational problem
$$
\inf_g \int |\hA|^2dx = \inf_g \int \sum_{\ell=1}^d\big|g^{-1}   A _\ell   g + g^{-1}   \p _\ell   g\big|^2  dx.
$$
It is easy to verify that the infimum is achieved and that $\sum^d_{\ell=1} \p _\ell   \hA _\ell =0$ \cite{DELLZ}. Thus the frame $\{\he, J\he \}$ is a Coulomb frame with $\he =g\cdot e$.

\subsubsection*{Remark} If $\phi \in W^{1,d}(\RR^d, M)$ and $M$ is compact then by a result of Schoen and Uhlenbeck \cite{SU83a, SU83b} $\phi $ can be approximated by smooth functions. Therefore by proposition \ref{global-f}  and  equation \eqref{exp}, the exponential  frame $\{e, Je\}$ on $\phi^{-1}TM$  belongs to $e\in W_{\rm loc }^{1,\frac d2}$. The local Coulomb gauge in \cite{Uh83} which satisfies
$$\begin{aligned}
&\p_\ell   A_k -\p_k    A_\ell
+ \lf[A_\ell, A_k\rt]
= F_{\ell  k}\in L^{\frac d2}\\
&\p_kA_k=0
\end{aligned}
$$
belongs to $L^d_x$ whence $e\in W_{\rm loc }^{1, d}$ for $d>2$.  For $d=2$ we need to require $\phi \in W^{1,p}$ for some $p>2$.

\subsubsection*{\bfi{GNLS}} The relation of Schr\"odinger maps to Gauge invariant Schr\"odinger equations is given through the frame coordinates
\begin{equation}\label{coor}\begin{cases}
&\partial_\alpha u = q_\alpha \cdot e,\\
&D_\alpha e = A_\alpha \cdot e.
\end{cases}
\end{equation}
for $\alpha = 0,1,\cdots, d$. Given such $\{q_\alpha, A_\alpha\}$,  let $   F_{ \alpha\beta} = F(q_\beta, q_\alpha) =R(q_\alpha\cdot e, q_\beta\cdot e)$ where $R$ denotes the Riemann curvature tensor of  $M$. We have,

\begin{proposition} Given a smooth Schr\"odinger map $u:\RR^d \times \RR \to M$ and a frame $\{e, Je\}$  on $u^{-1}TM$; the coordinates  $( q_\alpha, A_\alpha)$ for $\alpha = 0,1,\cdots, d$, given by \eqref{coor} satisfy
\begin{equation}\tag{GNLS}
\begin{cases}
&q_0 = i   \CD_\ell   q_\ell\\
&\CD_t  q_\ell = i  \CD^2_k   q_\ell + i   F_{\ell  k}  q_k,\\
&\CD_\ell   q_k= \CD_k  q_\ell, \\
&\p_\alpha   A_\beta -\p_\beta    A_\alpha + \lf[A_\alpha, A_\beta\rt] = F_{\alpha  \beta}
\end{cases}
\end{equation}
for $k, \ell = 1,\cdots, d$ and   $\alpha, \beta  = 0,1,\cdots, d$; and  where we summed on repeated indices.
\end{proposition}
\begin{proof}
Write $q=(q_0, q_1, \dots, q_d)\in \CC^{n\times (d+1)}$.  The $\CC^n$ valued functions
 $q_\alpha$, $\alpha=0,1, \dots, d$ have to satisfy
\begin{alignat}{2}
&q_0 = i   \CD_\ell   q_\ell &  \qquad \qquad&\qquad \qquad \text{since}\quad \partial_t u =J D_\ell \partial^\ell u \label{eq:sm}\\
&\CD_\alpha   q_\beta = \CD_\beta    q_\alpha& &\qquad \qquad\text{since }\quad
D_\alpha    \p_\beta   u= D_\beta  \p_\alpha   u
\label{eq:com}
\end{alignat}
The equations for the matrices $\{A_\alpha\}_{\alpha =0}^d$ can be derived from the curvature
equation
\begin{equation*}
\lf[D_\alpha, D_\beta\rt] e_a = R\big(\p_\beta    u, \p_\alpha   u\big)e_a = R\big(q_\beta\cdot e, q_\alpha \cdot e\big)e_a = F(q_ \alpha, q_\beta)\cdot e_a = F_{ \alpha\beta}\cdot e_a.
\end{equation*}
Note that $F_{\alpha  \beta}$ is  bilinear  in $(q_\alpha, q_\beta)$ and  is calculated from the Riemannian  curvature  and the frame on $u^{-1}TM$. 
Moreover in terms of the given frame we have
\begin{equation*}
\lf[\CD_\alpha, \CD_\beta\rt] = \p_\alpha   A_\beta -\p_\beta    A_\alpha
+ \lf[A_\alpha, A_\beta\rt]
= F_{\alpha  \beta}.
\end{equation*}
 Equations \eqref{eq:sm} and \eqref{eq:com} can be simplified by substituting \eqref{eq:sm} in equation \eqref{eq:com} for $\alpha = 0$ to obtain
$$
\CD_t  q_\ell = \CD_\ell  q_0 = i  \CD_\ell   D_k   q_k.
$$
By commuting $[\CD_\ell, \CD_k]$ and using the fact that $\CD_\ell   q_k= \CD_k  q_\ell$
we obtain the (GNLS) system.
\end{proof}
\subsubsection*{Remarks}
1. Given $u$ a solution of (SM) and a  choice of frames $\{e, Je\}$ we can compute $A_\alpha$ from $D_\alpha e_a= A_{a \alpha}^b e_b$.  By choosing another frame $\{\he , J \he\}$, where $\he_a = g^b_ae_b$ and  $g\in SU(n)$,  the connection  $\p_\alpha \he_a = \hA^b_{a \alpha} \he_b$ can be determined from $A_\alpha$ by
$
\p_\alpha g^b_a + A^b_{c\alpha} g^c_a = g^b_c\hA^c_{a\alpha}
$, or in matrix notation
$$
\hA_\alpha =g^{-1}A_\alpha g + g^{-1}\p_\alpha g\quad\quad \alpha = 0,1,\cdots, d.
$$
Thus  the equations for  $A_{\alpha }$ in  (GNLS) are  underdetermined unless we fix a choice of  the orthonormal basis $\{e, Je \}$.   Throughout this paper  we fix the frame by choosing  the Coulomb gauge which is given by $\sum_1^d\p_\ell A_\ell = 0$. \\

2. For $M$ a Riemann surface, the gauge group is $\UU(1)$, $q_\ell\in \CC$, $A_\alpha= ia_\alpha$ and $F(q_\alpha, q_\beta)= F_{\alpha \beta} = if_{\alpha \beta}$  where $a_\alpha, f_{\alpha \beta}\in \RR$.   In this case the (GNLS) system simplifies to
\begin{equation}\label{rs}\begin{aligned}
&\CD_tq_\ell= i\CD_k^2q_\ell - \kappa(u) \langle q_ \ell, i q_k\rangle q_k\\
&\CD_\ell q_k = \CD_k q_\ell\\
&\p_\ell a_j -\p_ja_\ell = f_{\ell j} = \kappa(u)\langle q_ \ell, i q_j\rangle\\
&\p_\ell a_0 -\p_ta_\ell =  f_{\ell 0} = -\kappa(u)\langle q_\ell ,  \CD_jq_j\rangle
\end{aligned}
\end{equation}
where $\kappa$ is the Gauss curvature of $M$, and where for two complex numbers $z$ and $w$ we used the notation $\langle z,w \rangle = \text{Re}(z\bar{w})$.  In this case it is always possible to put the above system in the Coulomb gauge globally by the gauge transformation $\hq_\ell =( \exp{i\theta})q$ and $\ha_\alpha = a_\alpha +\p_\alpha\theta$ where $\Delta \theta = -\p_\ell a_\ell$. In this Coulomb gauge equations \eqref{rs} transform into
\begin{equation}\label{rsc}
\begin{aligned}
&\CD_tq_\ell= i\CD_k^2q_\ell -\kappa(u) \langle q_ \ell, i q_k\rangle q_k\\
&\CD_\ell q_k = \CD_k q_\ell\\
&\Delta a_j = \p_k f_{kj} = \p_k (\kappa(u)\langle q_ k, i q_j\rangle)\\
&\Delta a_0 = \p_k f_{k0}=-\p_k\lf(\kappa(u)(\p_j \langle q_k ,  q_j\rangle - \half\p_k |q_j |^2)\rt)\\
&\p_k a_k = 0.
\end{aligned}
\end{equation}
3.  In general  the system (GNLS) depends on $u$ which appears in $F_{\alpha \beta}$.  For constant curvature $M$, the Schr\"odinger map $u$ does not appear explicitly in (GNLS). Thus we can consider the system (GNLS) on its own as an evolution problem. In this case the equations $\CD_\ell   q_k= \CD_k  q_\ell$ should be viewed as a compatibility condition which will be  satisfied under the evolutions of $q_\ell$ provided they are satisfied initially. Thus one of the questions we are interested in here is : Given $(q_\ell, A_\alpha)$  solutions of the (GNLS) in the Coulomb gauge, is there a Schr\"odinger map $u$  and a frame  $ \{e,J e \}$ such that  $  \p_\ell u = q_\ell \cdot e$ and
 $D_\alpha e = A_\alpha \cdot e$?

\section{Maps into $\SS^2$}\label{sphere}

\subsection*{\bfi{One space dimension}}A Schr\"odinger map $  u: \RR \times \RR \to \SS^2 \subset \RR^3 $ is a solution to
\begin{equation}\label{ll1}\begin{cases}
 & \partial_t u  =    J    D_{x} \partial_{x} u = u \times u_{xx}= \partial_x ( u \times u_x)    ,   \qquad u \in \SS^2\subset \RR^3\\
 &u(0) =u_0.
 \end{cases}
 \end{equation}
In this case the associated  GNLS system in the  Coulomb gauge $A_1=0$ is the nonlinear Schr\"odinger equation
$$
i \p_t q - \p^2_x q - \half   |q|^2   q= 0,
$$
and the transformation between $u$ and $q$ is given by
\begin{equation}\label{FS}
\begin{aligned}
&\p_t  u = p\cdot e = p_1  e+ p_2  u\times e \quad \qquad && u(0,x)
= u_0(x) \\
&\p_x u = q\cdot e = q_1  e+ q_2  u\times e\\
&D_t e = \p_t  e + p_1 u = - \frac12{|q|^2} u\times e && e(0,x) =e_0(x)\\
&D_x e =  \p_x  e + q_1 u = 0
\end{aligned}
\end{equation}
where $p=iq_x$.

For smooth solutions one can easily show the equivalence between solutions to the (SM) and solutions to the (NLS).
\begin{proposition} \label{pro:s}
1.  Given a smooth smooth solution $u$  of \eqref{ll1} there exist a frame  $\{e, u\times e\}$ for  $u^{-1}T\SS^2$   and a solution to the NLS \, \, $\p_t q= i\big( \p^2_x q + \half   |q|^2   q\big)$ such that
\begin{align*}
&\p_x u =q\cdot e = q_1e + q_2 u\times e \\
&D_x e  = \p_xe  + q_1u =0.\\
\end{align*}
\noindent 2. Conversely given a smooth solution $q$ to (NLS) with data $q_0$, a point $m \in \SS^2$ and $v_0 \in T_m\SS^2$ with $|v_0|=1$, there exists a unique solution $u$ to \eqref{ll1} and a frame $\{e, Je\}$ for  $u^{-1}T\SS^2$  such that \eqref{FS} is satisfied
with $u(0,0) = m$ and $e(0,0)=v_0$.
\end{proposition}

\begin{proof} 1. Let $u$ be a solution of \eqref{ll1},  $\{e, u\times e\}$ be any frame on $u^{-1}T\SS^2$, and let   $\p_\alpha u = q_\alpha\cdot e$ and $D_\alpha e =( i a_\alpha)\cdot e$.  Apply the gauge transformation $\p_x\theta = -  a_1$ to put the system in the Coulomb gauge $\ha_1 =0$. Since in this case the scalar curvature $\kappa=1$ we conclude from \eqref{rsc} that $\ha_0 = -\half|q|^2$ and that  $\hq$ satisfies
$$
(\p_t + i\ha_0)\hq = i\p_x^2 \hq
$$
which is the (NLS).

2. We start by constructing $u(0,x)
 \stackrel{def}{=} u_0(x)$ and $e(0,x)
\stackrel{def}{=} e_0(x)$ by solving the ODEs
\begin{align*}
&\p_x  u_0 = q_{01}  e_0 +  q_{02}  u_0 \times e_0\\
&\p_x   e_0 + q_{01}   u =0\\
&u_0 (0) =m \qquad e_0 (0) =v_0.
\end{align*}
where $q_{01} (x) + iq_{02} (x) = q_0(x)= q(0,x)
$. It is easy to check that $e_0(x) \perp u_0(x)$ and that
$|u_0 (x)|= |e_0(x)|=1$.

To construct $u$ and $e$ we evolve the data in time using  \eqref{FS}
\begin{equation*}
\begin{aligned}
&\p_t  u = p_1  e+ p_2  u\times e\\
&\p_t  e + p_1 u =- \frac12{|q|^2} u\times e\\
&u(0,x)
 = u_0(x) \qquad e(0,x)
=e_0(x).
\end{aligned}
\end{equation*}
where $p=iq_x$, to find  $u(t,x)
 \in \SS^2$ and $e \in T_u   \SS^2$,  $ |e(t,x)
|=1$.  To verify that $u$ solves $(SM)$ and that $D_x e =0$ we set
$\p_x  u= \tq \cdot e$ and $D_x e= a   u\times e$.
Then $\tq (0,x)
=q(0,x)
$ and $a(0,x)
=0$ by construction. By commuting
derivatives, we have
\begin{align*}
&D_t  \p_x u = D_x  \p_t  u \Rightarrow  \tq_t = i \p^2_x q
-q_xa +i   \frac 12{|q|^2}   \tq\\
& D_t   D_x   e - D_x   D_t  e = R(u_t, u_x)  e \Rightarrow
\p_t   a+ \p_x   \frac 12{|q|^2}  =\tq_1   q_{1x} + \tq_2  q_{2x}.
\end{align*}
Therefore $\tq-q$ and $a$ satisfy the ODEs
\begin{align*}
&\p_t(\tq -q) =i  \frac{|q|^2}{2} \big(\tq-q\big) + ipa\\
&\p_t  a = \big(\tq_1 - q_1 \big) q_{1x} + \big(\tq_2 - q_2\big)q_{2x},\\
&\big(\tq-q\big) (0,x)
=0 \qquad a(0,x)
=0,
\end{align*}
which imply  $\tq \equiv q$ and $a\equiv 0$.  Since $\p_t u =(iq_x)\cdot e$ and $D_x e =0$  we conclude that
$u$ solves  \eqref{ll1}.  The uniqueness of $u$ follows from the uniqueness
of the solutions to the ODEs and NLS.
\end{proof}
For rough initial data we can show equivalence of solutions under weak integrability conditions.
\begin{theorem}\label{thm:rgh}
Let
$
q\in L^2_{|t|<T}L^2_{x, \rm loc}$ satisfying
$|q|^2 \in  L^2_{|t|<T}(  H^{-1})
$
be the limit of smooth solutions, i.e.,  $\exists
q_k$ smooth solutions of (NLS) such that
$$
q_k \to q \in L^2_{tx, \rm loc} \quad\text{and}\quad
|q_k|^2 \to |q|^2 \in L^2_{t, \rm loc }  H^{-1}.
$$
Then there exists a weak solution $u\in L^2_{|t|<T}( H^1_{\rm loc})\cap C_{|t|<T}(L^2_{\rm loc})$
to \eqref{ll1} and a frame $\{e, u\times e\}$ of $u^{-1}  T\SS^2$
such that $e\in  L^2_{|t|<T}( H^1_{\rm loc})\cap C_{|t|<T}(L^2_{\rm loc})$.  Moreover the solution is unique up to isometries on the sphere.
\end{theorem}
\subsection*{Remarks}  1)  In one dimension, Vargas-Vega \cite{VV} showed local well posedness for the cubic NLS in a space containing $L^2$ and scaling like ${\dot H}^{-1/3}(\RR)$. Their solutions belong to $L^3_{|t|< T} L^6_x$ ( or $L^4_{|t| < T} L^4_x$). The critical scaling for the 1d NLS is that of ${\dot H}^{-1/2}(\RR)$.   Below $L^2$ however, the Galilean transformations are not preserved and the problem is ill-posed in the Sobolev class \cite{KPV}.

2)  In \cite{GRV} it is proved that a vortex filament can develop a singularity in the shape of a corner in finite time. This shows the existence of a Dirac delta singularity for the corresponding cubic NLS solution. For NLS data in $L^2$ such a behavior cannot occur due to  mass conservation.

\medskip

\begin{proof}
By proposition \ref{pro:s} we can construct from $\{q_k\}$ smooth solutions  $u_k$   of
\eqref{ll1} and frames $\{e_k, u_k \times e_k\}$ of $u^{-1}_k   T\SS^2$
such that
\begin{equation}\label{FSK}
\begin{aligned}
&\p_t  u_k = p_k\cdot e_k = p_{k1}  e_k+ p_{k2} \, u_k\times e_k \quad \qquad && u_k(0,x)
= u_0(x)\\
&\p_x u_k = q_k\cdot e_k = q_{k1}  e_k+ q_{k2} \, u_k\times e_k\\
&D_t e_k = \p_t  e_k + p_{k1} u_k = - \frac12{|q_k|^2}\, u_k\times e_k && e_k(0,x)
=e_0(x)\\
&D_x e_k =  \p_x  e_k + q_{k1} u_k = 0
\end{aligned}
\end{equation}
where $p_k=i q_{kx}$.  By the hypothesis of the theorem we can pass to the
limit in \eqref{FSK} and thus $u$ and $e$ satisfy equation \eqref{FS} in
the sense of distributions.  From the equations for $\p_x u$ and $\p_x e$ we conclude that $u$ and $e$ are in
$ L^2_{|t|<T}( {H}^1_{\rm loc})$.   From the equations for $\p_t u$ and $\p_t e$ we have  $\forall \varphi \in C^\infty_0(\RR)$ $\varphi u$ and $\varphi e$ are in $C_{|t|<T}(H^1)$.  From computing
$$
\p_t \int |u(t,x)-u(t_0,x)|^2\varphi(x)dx =2 \int \langle p\cdot e , u(t,x)-u(t_0,x)\rangle \varphi(x)dx,
$$ where $\langle \, , \, \rangle$ is the inner product in $\Bbb R^3$,
and a similar expression for $e$ we conclude that $u$ and $e$ are in $C_{|t|<T}(L^2_{\rm loc})$.  Note that in this case \eqref{FS} implies that for every $t\in(-T,T)$,  $u(t,\cdot)$ and $e(t,\cdot)$ are in $\dot H^1$.

To show uniqueness, let $(u,e)$ and $(\tu, \te)$ be two such solutions that  satisfy \eqref{FS}. Then by using the isometries on $\SS^2$ we can assume that $u(0,0)=\tu(0,0)$ and $e(0,0)=\te(0,0)$. Equation \eqref{FS} implies that
$u(0,x)=\tu(0,x)$ and $e(0,x)=\te(0,x)$ for all $x$.   Set $\delta u = u-\tu$,
$\delta e = e - \te$, and $\delta f = u \times e-\tu \times \te$ then
\begin{equation}\label{FSD}
\begin{aligned}
\p_t   \delta u &= p_1  \delta e + p_2   \delta f\\
\p_t   \delta e &= -p_1  \delta u - \half   |q|^2  \delta f\\
\p_t   \delta f &= -p_2  \delta u + \half   |q|^2  \delta e\\
\delta u(0) &=0, \qquad \delta e (0,x)
 = \delta f (0,x)
 =0
\end{aligned}
\end{equation}
which can be written in matrix notation as
$$
\p_t V = BV \qquad V(0)=0.
$$
Since $B$ is skew symmetric and is locally in $L^2  H^{-1}$ and $V$ is
locally in $L^2  H^1 \cap L^\infty$ then for any $C^\infty_0 \ni \varphi$
$$
\frac d{dt}\int |V(t,x)|^2\varphi(x)dx = 2 \int  \langle BV, V\rangle \varphi  dx =0
$$
and therefore $V\equiv 0$.
\end{proof}

\subsubsection*{ \bfi{ Higher dimensional maps into $\mathbf{\SS^2}$}}
A Schr\"odinger map $  u: \RR^d \times \RR \to \SS^2 \hookrightarrow \RR^3 $ is a solution to
\begin{equation}\label{ll}
\partial_t u  =     u \times \Delta u = \partial_{x_j} ( u \times \partial_{x_j} u ).
\end{equation}
In this case, since $n=1$, we have $q_\alpha\in \CC$, $A_\alpha= ia_\alpha$ and $F(q_\alpha, q_\beta)= F_{\alpha \beta} = if_{\alpha \beta}$  where $a_\alpha, f_{\alpha \beta}\in \RR$.   Given a Schr\"{o}dinger map $u$ into $\SS^2$ and a frame $\{e,Je\}$ we recall  (GNLS)  for $(q_k,a_k)$
\begin{align*}
&\CD_tq_\ell= i\CD_k^2q_\ell  -  i\langle q_ \ell, i q_k\rangle q_k\\
&\CD_k q_\ell = \CD_\ell q_k\\
&\p_\ell a_j -\p_ja_\ell = f_{\ell j} = \langle q_ \ell, i q_j\rangle\\
&\p_\ell a_0 -\p_ta_\ell =  f_{\ell 0} = -\langle q_\ell ,  \CD_jq_j\rangle
\end{align*}
and the transformation between $u$ and $q$
\begin{equation}\label{gll}\begin{aligned}
&\p_t u = q_0\cdot e\\
&\p_{\ell} u = q_\ell\cdot e\\
&D_t e = a_0 u\times e\\
&D_\ell e = a_\ell u\times e
\end{aligned}
\end{equation}
where $q_0= i\CD_kq_k$.  In the Coulomb frame this system simplifies to
\begin{equation}\label{gllc}
\begin{aligned}
&\CD_tq_\ell= i\CD_k^2q_\ell -i \langle q_ \ell, i q_k\rangle q_k\\
&\Delta a_j = \p_k f_{kj} = \p_k\langle q_ k, i q_j\rangle\\
&\Delta a_0 = \p_k f_{k0}= -\p_\ell\p_j \langle q_\ell ,  q_j\rangle + \half\Delta |q_j |^2
\end{aligned}
\end{equation}
along with the Coulomb frame equation and compatibility conditions
\begin{equation} \label{E:compat} \begin{aligned}
&\p_k a_k=0\\
&\CD_k q_\ell = \CD_\ell q_k\\
&\p_\ell a_j -\p_ja_\ell = f_{\ell j} = \langle q_ \ell, i q_j\rangle\\
&\p_\ell a_0 -\p_ta_\ell =  f_{\ell 0} = -\langle q_\ell ,  \CD_jq_j\rangle.
\end{aligned}\end{equation}
It is easy to verify that \eqref{E:compat} are satisfied by smooth solutions of~\eqref{gllc} for all $t$ if they are satisfied at $t=0$ and $a$ decays at infinity.

\begin{proposition}\label{eqsd}
Given a smooth solution $(q, a)$ of \eqref{gllc} satisfying \eqref{E:compat}, a point $m\in \SS^2$
and a vector $v_0 \in T_m  \SS^2$ with $|v_0|=1$, then there exists a
unique solution $u$  of the \eqref{ll}  system and a frame $\{e, Je\}$ for $u^{-1}  T\SS^2$ such that \eqref{coor}, i.e.,
\begin{equation*}
\begin{aligned}
\p _\alpha u&= q _\alpha\cdot e\\
D _\alpha  e &= a _\alpha u\times  e
\end{aligned}
\end{equation*}
holds with $u(0,0)=m$ and $e(0,0) = v_0$.
\end{proposition}

\begin{proof}
Given  $(q, a)$ solution to \eqref{gllc} we first construct the initial data
for $u$ and for the frame $\{e, u\times e\}$.  This will be done inductively on every coordinate $x_1, x_2, \cdots, x_d$.  We start by solving
\begin{align*}
&\p_1  w_1 (x_1)  = q_1(x_1, 0,\dots, 0\big)
\cdot e_1(x_1) \\
&D_1  e_1 (x_1)  = \p_1  e_1 (x_1)  + \langle e_1 (x_1) , \p_1  w_1 (x_1) \rangle  w_1(x_1)
= a_1(x_1, 0,\dots, 0\big) w_1\times e_1 (x_1) \\
&w_1(0)=m,\qquad e_1(0)= v_0.
\end{align*}
It is easy to verify that $\{e_1, w_1\times e_1\}$ are a frame
along the curve $w_1^{-1}T\SS^2$.  Repeat this process to construct $w_2(x_1,x_2)$ and $\{e_2, w_2\times e_2\}$ from
\begin{align*}
&\p_2  w_2   = q_2(x_1, x_2, 0,\dots, 0\big)
\cdot e_2\\
&D_2  e_2 = \p_2 e_2 + \langle e_2 , \p_1  w_2\rangle  w_2
= a_2(x_1,x_2, 0,\dots, 0\big) w_2\times e_2\\
&w_2(x_1, 0)=w_1(x_1),\qquad e_2(x_1,0)= e_1(x_1).
\end{align*}
This construction terminates by constructing $u_0\sdef w_d(x_1,\cdots, x_d)$ and $\{e_0, u_0\times e_0\}\sdef \{e_d, w_d\times e_d\}$.

To verify that  $u_0$ and $\{e_0, Je_0\}$ satisfy \eqref{coor} at $t=0$, we note  that by construction the equations hold on $\big(x_1, 0, \dots,0 \big)\in \RR^2$. Here $\RR^2$ denotes the $x_1x_2-$plane. To show that the same holds for on $\big(x_1, x_2, 0\dots ,0\big)\in \RR^2$ we compute
\begin{equation}\label{A}
\begin{aligned}
D_2  \p_1 u_0 &= D_1  \p_2  u_0\\
D_2  D_1 e_0 &= D_1  D_2  e_0 + R\big(\p_1 u_0,\p_2 u_0\big)e_0.
\end{aligned}
\end{equation}
By our construction we have
\begin{gather}\label{B}
\begin{aligned}
&\p_2 u_0 = q_2 \cdot e_0\quad &D_2 e_0=a_2 \, u_0 \times e_0\quad & (x_1, x_2) \in \RR^2\\
&\p_1 u_0 = \hq_1 \cdot e_0\quad &D_1 e_0=\ha_1 \, u_0 \times e_0\quad & (x_1, x_2) \in \RR^2
\end{aligned}\\
\hq_1(x_1, 0) = q_1(x_1,\cdots,0) \qquad  \ha_1(x_1, 0)=a_1(x_1, \cdots, 0) \notag
\end{gather}
Substituting \eqref{B} in \eqref{A} we obtain the following ODEs
\begin{align*}
\p_2 \big(\hq_1 - q_1\big) + ia_2 \big(\hq_1 -q_1\big) = i\big(\ha_1 - a_1\big)q_2\\
\p_2 \big(\ha_1 - a_1\big)  = \hf_{12} - f_{12} = \langle \hq_1 - q_1, q_2\rangle
\end{align*}
where $\hf_{12}  = \big(\hq_1, iq_2\big)$.  Since $0$ is a solution of this
ODE by uniqueness we have $\hq_1=q_1$ and $\hq_2 =q_2$ in the $x_1x_2$-plane. Repeating this process for $x_3,\cdots, x_d$, we obtain the desired result.

To construct $u(t,x)
$ and $\{e, u\times e\}$ we solve the ODEs
\begin{equation*}\begin{aligned}
&\p_t u = q_0\cdot e\\
&D_t e =\p_te +q_{01}u=a_0 u\times e\\
&u(0,x)
=u_0(x)\qquad e(0,x)
=e_0(x)
\end{aligned}
\end{equation*}
To verify that $u$ solves $(SM)$ and \eqref{coor} holds we set
$\p_\ell  u= \tq_\ell \cdot e$ and $D_\ell e= \ta_\ell   u\times e$ and define $\tilde{\CD}_\ell = \p_\ell + i \ta_\ell$. Then $\tq (0,x) =q(0,x) $ and $\ta(0,x)=a(0,x)$ by construction. By commuting
derivatives, we have
\begin{align*}
&D_t  \p_\ell u = D_\ell  \p_t  u \Rightarrow \CD_t \tq_\ell= i \tilde{\CD}_\ell q_0  = i\tilde{\CD}_\ell \CD_k q_k\\
& D_t   D_\ell    e - D_\ell    D_t  e = R(u_\ell, u_t )  e \Rightarrow
\p_t   \ta_\ell - \p_\ell   a_0 =\tilde{f}_{0\ell},
\end{align*}
where  $\tilde{f}_{0\ell} = -i F(q_0, \tq_\ell)$.
Therefore $\tq_\ell-q_\ell$ and $\ta_\ell -a_\ell$ satisfy the ODEs
\begin{align*}
&\CD_t(\tq_\ell-q_\ell) =-(\ta_\ell -a_\ell)\CD_k q_k\\
&\p_t ( \ta_\ell -a_\ell)  = \tilde{f}_{0\ell}  -f_{0\ell} = \langle q_0, i(\tq_\ell -q_\ell) \rangle,\\
&(\tq_\ell-q_\ell)(0,x) =0 \qquad (\ta_\ell -a_\ell)a(0,x) =0.
\end{align*}
which imply  $\tq_\ell \equiv q_\ell $ and $\ta_\ell \equiv a_\ell $, and thus   we conclude that
$u$ solves  \eqref{ll1}.  The uniqueness of $u$ follows from the uniqueness
of the solutions to the ODEs and NLS.
\end{proof}

\begin{theorem}
Given a solution $q_k$ to \eqref{gllc} such that
\[ \begin{split}
q_k \in C ( [0,T], L^2(\RR^d)) \cap L^6 ( [0,T], L^3(\RR^d))
&\quad \text{for}\ d= 2\\
q_k \in C ( [0,T], L^2(\RR^d) \cap L^d (\RR^d) )
& \quad \text{for}\ d \geq 3
\end{split}
\]
and assume that $q_k$ is the $C( L^2)$ limit of smooth solutions,
i.e., $\exists \{ q_{k_{(j)}} \} $ where $q_{k_{(j)}} \to q_k$ in $C\big([0,T],
L^2_x\big)$.
Then  there exists  a solution $u \in L^\infty( \dot H^1)\cap C(L^2_{\rm loc})$ of the Schr\"odinger maps equation \eqref{ll} and a frame $\{e,Je\}$, where $e  \in L^\infty( \dot H^1)\cap C(L^2_{\rm loc})$ for $d>2$ and $e \in L^\infty (\dot H^1_{\rm loc}) \cap C(L^2_{\rm loc})$ for $d=2$, that
satisfies the coordinates equation \eqref{coor}
$$\begin{cases} \partial_{\alpha} u   =    q_{\alpha}\cdot e   \\
D_\alpha e = \partial_{\alpha} e    +     q_{\alpha 1} u   =    a_{\alpha} u\times e,
 \end{cases} $$
Moreover if there are two such $\{u, e\}$ and $\{ \tilde{u}, \tilde{e}\}$
that have the same coordinates given by  \eqref{coor} with initial data $u(0, x_0) = \tilde{u} (0, x_0)$,
$e(0, x_0)= \tilde{e} (0, x_0)$ at one Lebesgue point $(0, x_0)$,
$x_0 \in \RR^d$ of $u -\tilde u$ and $e - \tilde e$, then $u \equiv \tilde{u}$ and $e \equiv \tilde{e}$.
\end{theorem}

\subsection*{Remarks} 1) The assumption that $q_k \in C (L^2) \cap L^6 (L^3)$ in
two dimensions guarantees finite energy plus a Strichartz norm.
This is necessary to make sense of all the terms in \eqref{coor},
such as $a_0 $ and is not needed for the existence of weak solutions . Other Strichartz choices are also possible.\\
2) The assumption that $q_k \in C (L^d)$ for $d \geq 3$ is
much weaker than the space $C (H^{\frac d 2 -1})$
which is the optimal space for existence of solutions to \eqref{gllc}.\\
3) The assumption $\tilde u = u$ and $\tilde e = e$ at a Lebesgue point in the uniqueness statement can also be replaced by any decay to $0$ of $\tilde u(0, \cdot) -u (0, \cdot)$ and $\tilde e(0, \cdot) - e(0, \cdot)$ as $|x|\to \infty$.

\begin{proof} From the expression for $a_0$ and $a_j$ in \eqref{gllc}
$$\begin{aligned}
&\Delta a_j = \p_k f_{kj} = \p_k\langle q_ k, i q_j\rangle\\
&\Delta a_0 = \p_k f_{k0}= -\p_\ell\p_j \langle q_\ell , i q_j\rangle + \half\Delta |q_j |^2
\end{aligned}
$$
we have $a_j\in C(L^d)$,  $a_0 \in C(L^{d/2})$ for $d > 2$  and $a_j\in L^3(L^{6})$,  $a_0\in L^3(L^{3/2})$ for $d=2$.
By proposition \ref{eqsd} we can construct from $\{q_{k(j)}\}$ smooth solutions  $u_{(j)}$   of
\eqref{ll} and frames $\{e_{(j)}, u_{(j)} \times e_{(j)}\}$ of $u^{-1}_{(j)}   T\SS^2$
such that
\begin{equation}\label{FSKd}
\begin{aligned}
\p_\alpha   u_{(j)} &=q_{\alpha (j)} \cdot e_{(j)}\\
\p_\alpha  e_{(j)} &= -q_{\alpha1(j)}   u_{(j)}  +a_{\alpha (j)}   u_{(j)} \times e_{(j)}
\end{aligned}
\end{equation}
By the regularity hypothesis on $q_k$ given in  the theorem we can pass to the
limit in \eqref{FSKd} and thus $u$ and $e$ are in $L^\infty(\dot{H}^1)\cap C(L^2_{\rm loc})$ (if $d=2$, $e\in L^\infty(\dot H^1_{\rm loc})$), they  satisfy equation \eqref{coor} in
the sense of distribution, and $u$ solves \eqref{ll}.

To show uniqueness assume $u$ and $\tilde u$ are
two solutions that satisfy \eqref{coor} and agree at a point say $(0,0)$.
We first show that the data for $u$ and  $\tilde u$
are the same.  Let $f= u \times e$ and $\tilde f = \tilde u\times  \tilde e$,
then from \eqref{coor} we have at $t =0$
\[ \begin{split}
\partial_k u & = q_k \cdot e = q_{k1} e + q_{k2} f\\
\partial_k e &= -q_{k1} u + a_k f\\
\partial _k f &= -q_{k2} u - a_k f
\end{split} \]
and the same for $\{ \tilde u, \tilde e, \tilde f \}$. Since $q_k(0, \cdot) \in L^d (\RR^d)$, for $d=2$, it is straight forward to prove that $a_k (0, \cdot)$ is in the dual space of $L^r \cap L^{\frac {2r}{r+2}}$ for any $r>2$ and for $d>2$, $a_k(0, \cdot) \in L^d$. Thus $a_k (0, \cdot)$ is in $L^1_{\rm loc}$. Therefore we may take differences in the above linear equations to obtain
\[
|u - \tilde u|^2 + |e - \tilde e|^2 + |f - \tilde f|^2 = \text{constant}.
\]
Since $u(0, 0) = \tilde u (0, 0)$ and
$e(0, 0) = \tilde e (0, 0)$ then
$u \equiv \tilde u$ and $e \equiv \tilde e$ at $t=0$.

To show that $u \equiv \tilde u$ for all $t$ we use the time
derivative part of \eqref{coor}
\[ \begin{split}
\partial_t u &= p \cdot e = q_{01} e + q_{02} f\\
\partial_t e &= - q_{01} u + a_0 f\\
\partial_t f &= -q_{02} u - a_0 e
\end{split}  \]
and the same for $\{ \tilde u, \tilde e, \tilde f \}$.
Again since $p$ and $a_0$ are in  $L^1_{\rm loc}(H^{-1})$ we have
\[
 \partial_t \big( |u - \tilde u|^2 + |e - \tilde e|^2 +| f - \tilde f|^2 \big) = 0
\]
and since at $t=0$, $u = \tilde u$ and $e = \tilde e$,
then $u \equiv \tilde u$ and $e \equiv \tilde e,  \forall (t,x) \in \R \times \RR^d$.
\end{proof}

\section{  \sms\ into $\HH^2$ } \label{hyperbolic}

The Cauchy problem for Schr\"odinger maps into the hyperbolic plane
\( u : \RR^d \times \RR \to \mathbb{H}^2 \)
has two difficulties that are not present when the target is $\mathbb{S}^2$.
The first difficulty  is due to the fact that $\mathbb{H}^2$ cannot be embedded isometrically and equivariantly in $\RR^k$.  The second is due to the non compactness of $\mathbb{H}^2$,
which makes controlling $u$ an issue.

The  first difficulty can be avoided by embedding $\mathbb{H}^2$ in the Lorentz space
\(( \RR^3, \eta )\) where\linebreak[4]
 $\eta = {\rm dia} ( -1, 1, 1)$ and the embedding is given
by
\[
\mathbb{H}^2 = \big\{ u \in (\RR^3 , \eta ) ; - u^2_0 + u^2_1 + u^2_2 = -1,\, u_0 >0
\big\} .
\]
The embedding is isometric and equivariant as becomes apparent
after introducing the coordinates $u_0 = \cosh \chi$,
$u_1 = \sinh \chi \cos \theta $
and $u_2 = \sinh \chi \sin \theta$.  The tangent space and the normal space for this embedding are given by
\begin{align*}
T_u \mathbb{H}^2 & = \{ v \in \RR^3 ; \< \eta u, v \> =0 \}\\*
N_u \mathbb{H}^2 &= \{ \gamma u ; \gamma \in \RR \} .
\end{align*}
The unit  normal at $u \in \mathbb{H}^2$ is the vector $u$ since
\( \< \eta u, u \> = -1 \).
For a vector $v \in T_u \h^2$ we introduce the notation
\[
\|v\|^2 = | \< v ,  \eta v \>|  = - v_0^2 +v_1^2 +v_2^2 ,
\]
and for a map $u : \RR^d \to \h^2 $
with $w \in u^{-1} T \h^2$
\[
\| w\|_{L^2}^2 = \int \| w(x) \|^2 \, dx.
\]

Given a map $\phi : \RR^d \to \HH^2 \subset \RR^3$,
the covariant derivative on $\phi^{-1} T \mathbb{H}^2$ is given by
\[
D_k V = \partial_k V - \< V, \eta \partial_k \phi \> \phi .
\]
The complex structure on $ T \mathbb{H}^2$ can be represented by
\[
J v = \eta (u \times v )
\]
where $\times$ is the usual cross product on $\RR^3$.
This is a consequence of
\( \< u, \eta J v \> = \< v, \eta J v \> = 0 \)
and $J^2 = -I$.

Using  the embedding $\h^2\subset (\RR^3,\eta)$  \sms\
\( u : \RR^d \times \RR \to \mathbb{H}^2 \)
can be written in divergence form as
\begin{equation} \label{E:smh}
\frac{\partial u}{\partial t} =
\eta ( u \times ( \Delta u - \< \nabla u , \eta \nabla u \> ) )
= \eta (u \times \Delta u)=  \eta \partial_k(u \times \partial_k u)
\end{equation}
or equivalently
\begin{equation}\label{smh}
\eta \bigg( u \times \frac{\partial u}{\partial t}\bigg)  =
- \Delta u + (\nabla u , \eta \nabla u) u
\end{equation}
In hyperbolic coordinates this system reduces to
\[ \begin{split}
(\sinh \chi ) \theta_t
 & = \Delta \chi - \sinh \chi \cosh \chi
| \nabla \theta |^2\\
(\sinh \chi ) \chi_t
 & = - \op{div} (\sinh^2 \chi \nabla \theta ) .
\end{split} \]

Given a smooth solution  to \eqref{smh} we can easily construct a frame $\{e\}$ in the Coulomb gauge and from section \ref{frame} the coordinates $\p_\ell u = q_\ell e$ satisfy
\begin{equation}\label{hypq}\begin{aligned}
&\CD_tq_\ell= i\CD_k^2q_\ell + i\langle q_ \ell, i q_k\rangle q_k\\
&\CD_\ell q_k = \CD_k q_\ell\\
&\p_\ell a_j -\p_ja_\ell = f_{\ell j} = -\langle q_ \ell, i q_j\rangle\\
&\p_\ell a_0 -\p_ta_\ell =  f_{\ell 0} = - \langle q_\ell , iq_0\rangle\\
&\p_k a_k=0
\end{aligned}
\end{equation}
where $q_0 =  i \CD_jq_j$. Conversely given a solution to \eqref{hypq} one can  repeat the construction given for the sphere in proposition \ref{pro:s} to obtain

\begin{proposition}\label{p2}
Given a smooth solution to \eqref{hypq},
a point $m \in \h^2$ and a vector $v_0 \in T_m \h^2$ with $\|v_0\|=1$,
then there exists a unique smooth solution to the \sms\ equation
\[
\begin{split}
\partial_t u & = \eta \partial _\ell (u \times \partial _\ell u )\\
u &\in \h^2 \subset (\RR^3, \eta )
\end{split}
\]
and a frame $\{e,Je\}$ for $u^{-1} T \h^2$ such that
$ u (0,0) = m$, $e (0,0) = v_0$, and \eqref{coor} holds.
\end{proposition}

\subsection*{\bfi{Weak finite energy solutions from $\RR^{2+1}$ into $\HH^2$}}

The difficulty of the non compactness of $\HH^2$ appears in constructing weak solutions and it can be overcome by requiring the map $u$ to converge to a point as $x\to \infty$. In particular, fix a point   $o\in \mathbb{H}^2$ and embed   $\mathbb{H}^2$ into Lorentz space with $o\to (1,0,0)$.  We will consider maps $u: \RR^2\to \h^2\subset (\RR^3,\eta)$  such that $u\to (1,0,0)$ as $x\to \infty$ and
$$
\int (u_0-1)dx =\int(\cosh\chi-1)dx <\infty.
$$
This is a reasonable assumption since, like the energy
\[
\|\nabla u\|_{L^2}^2 = \int |\nabla \chi|^2 + \sinh^2 \chi |\nabla \theta|^2 dx,
\]
$\int (u_0-1)dx$ is also a conserved quantity of the \sms.

Consider the Cauchy problem
\begin{equation}\label{smhc}
\begin{gathered}
\partial _t u = \eta ( u \times \Delta u )=  \eta \partial_k(u \times \partial_k u)\\
\|\nabla u(0)\| \in  L^2 (\RR^2 ), \quad u_0 (0) -1 \in L^1 (\RR^2 ).
\end{gathered}
\end{equation}
Since the equation is in divergence form then it is easy to conclude that the weak limit of finite energy smooth solutions is a weak  solution.

\begin{proposition}\label{lss}
Let $\{ u_k\}$ be a sequence of smooth solutions to the \sms\
equations \eqref{smhc} such that
\[
\|\nabla u_k\|_{L^2}^2\le C \qquad \text{ and }\qquad \int (u_{0_k} -1) dx \le C
\]
then $\exists$ a subsequence that converges w* to a weak solution of \eqref{smhc} $u \in L^\infty(W^{1,p}_{\rm loc})$ for any $p <2$.
\end{proposition}

\begin{proof}
From conservation of energy and the divergence form of the equation we have
\[
\int \|\nabla u_k (t)\|^2 dx \le C, \quad
\int ( u_{0k} (t) -1 ) dx \le C .
\]
In hyperbolic coordinates we have
\[
\int |\nabla \chi_k (t) |^2 \, dx \le C \quad
\int (\cosh \chi_k (t) -1) dx\le C .
\]
Thus \( | \chi_k (t) |_{H^1(\RR^2)} \le C \)
and from Moser-Trudinger inequality we have $\forall$ compact sets
$\Omega \subset \RR^2$
\[ \intl_\Omega e^{a\chi^2_k (t)} \,dx \le
C(\Omega, |\chi_k|_{H^1} ) \le C (\Omega ),
\]
for some $a>0$. These bounds on $\chi_k(t)$ imply the following Euclidean bounds on
$u_k(t)$
$$
|u_k(t)|_{L^p(\Omega)}\le C \quad \forall\ 1 \le p< \infty
$$
which in turn gives the Euclidean bounds
\begin{alignat*}{2}
&|u_k|_{ L^\infty (W^{1,p} (\Omega) ) } &\le  C\quad &\forall\ 1 \le p<2 \\
&|\partial_t u_k|_{ L^\infty (W^{-1,p} (\Omega) )} & \le  C\quad &\forall\ 1 \le p < 2  .
\end{alignat*}
Thus by going to a subsequence and a diagonalization argument we have
$\forall$ compact $\Omega \subset \RR^2$
\begin{alignat*}{2}
u_k \rightharpoonup u  &\in L^\infty (W^{1,p} (\Omega) ) &\quad & 1 \le p < 2 \\
u_k \to u &\in C(L^p (\Omega) ) & & 1 \le p < \infty .
\end{alignat*}
and this implies
\[
u_k \wedge \nabla u_k \rightharpoonup u \wedge \nabla u \in L^\infty
( L^p (\Omega) ) \quad 1 < p < 2.
\]
From the above and Fatou's lemma we conclude that $u$ is a weak solution of the \sms\ equation
with
\[
\int ( u_0 (t) -1 ) dx \le C.
\]
In order to show
\[
\int \|\nabla u (t)\|^2  dx \le C,
\]
we take an isometric embedding $\Phi: \HH^2 \to \RR^n$ satisfying $\Phi(o) =0$ and consider $\tilde u_k = \Phi \circ u_k$. Since $\chi$ is the geodesic distance to $o$ on $\HH^2$ and the intrinsic metric $\|\cdot \|$ on $T\HH^2$ coincides with the metric induced by $\Phi$, we have $|\tilde u_k (t) |_{H^1}^2 \le C$. Due to the pointwise convergence of $u_k$ to $u$, we have $\tilde u_k \rightharpoonup \tilde u \triangleq \Phi \circ u$ in $H^1$ and \[
\|\nabla u\|_{L^2}^2 = |\nabla \tilde u|_{L^2}^2 \le C.
\]
\end{proof}

To construct a sequence $\{u_k\}$ such that  $\partial_t u_k  -\eta \partial _\ell (u_k \times \partial _\ell u_k )\to 0$ in the sense of distribution and such that $\|\nabla u_k(t)\|_{L^2}< C$ and  $ |u_{0k} (t) -1|_{L^1}<C $, we introduce the parabolic perturbation
\[
\varepsilon \partial_t u - \eta ( u \times \partial_t u )  = \Delta u - \< \nabla u , \eta \nabla u \> u
\]
and show by using the frame coordinates $q$s that the above equation has global smooth solutions with the desired bounds.

\begin{proposition}\label{dgl}
Given $\varepsilon > 0$ and a function $u_* (y) \in \h^2$ such that
\[
u_* - (1, 0, 0 ) \in L^1 (\RR^2), \quad
\nabla u_* \in L^2 (\RR^2, T \h^2 ), \quad
D \partial u_* \in L^2 (\RR^2, T \h^2 ),
\]
there exists a unique global classical solution to
\begin{equation}\label{pert}
\begin{split}
\varepsilon \partial_t u - \eta ( u \times \partial_t u )
& = \Delta u - \< \nabla u , \eta \nabla u \> u \\
u (0, y ) & = u_* (y) \in \h^2,
\end{split}
\end{equation}
such that $u\in \h^2$ and
$$
\int\|\nabla u(t)\|^2 dx \le C, \quad \int ( u_0 (t) -1 ) dx \le C, \quad
\varepsilon \int\int\|\p_t u(t)\|^2 dx dt\le C, \quad\text{and}\  \varepsilon \int\int\|D \p  u(t)\|^2 dx dt \le C.
$$
\end{proposition}

\begin{proof} To show that solutions to equation \eqref{pert} stay in $\h^2$ we take the inner product of the equation with $\eta u$ to obtain
$$\begin{aligned}
&\frac 12 \varepsilon \p_t \< u, \eta u\> = \half \Delta  \< u, \eta u\>  +\< \nabla u , \eta \nabla u \>( 1+\< u, \eta u\>)\\*
&\< u, \eta u\>|_{t=0} = -1.
\end{aligned}
$$
which implies that  $u(t)\in \h^2$ .  To construct solutions
let $q$ be the Coulomb frame coordinates of $\p u$, then
\begin{equation}\label{hypc}
\begin{split}
(\varepsilon - i ) \CD_t q_\ell
& = \CD^2_k q_\ell -i \< q_\ell, i q_k \> q_k\\
\Delta a_j &= - \partial_k \< q_k, iq_j \>\\
 \Delta a_0 & = \p_\ell \langle q_\ell ,  q_0\rangle.
\end{split}
\end{equation}
where $(\varepsilon - i )q_0 = \CD_j q_j$.   By standard fixed point argument system \eqref{hypc} has local smooth solutions for initial data in $H^s$ for $s$ sufficiently large. Moreover the system has a conserved energy which can be  obtained by dividing the above equation by $(\varepsilon - i )$,
multiplying by $\bar q _\ell$ and taking the real part
\[
\frac d{dt} \int \frac 12 |q_\ell |^2
= \frac{- \varepsilon}{1 + \varepsilon^2}
 \int |\CD_k q_\ell|^2 + |\< q_k, iq_\ell \>|^2 .
\]
This implies global bounds
\[
E_0 = \frac 12 \int |q_\ell (t) |^2 \,dx
+ \frac \varepsilon{1+\varepsilon^2} \int^t_0 \int |\CD_k q_\ell (t) |^2
+ | \< q_k (t) , i q_\ell (t) \> |^2 \, dx \,dt .
\]

We will obtain the $H^1(\RR^2)$ estimate on $q$ by looking at $\CD q$. In fact,
\[
|\partial_k q_\ell |^2 \le |\CD_k q_\ell |^2  + |a_k q_\ell |^2 .
\]
Using the equation for $a_k$ and Sobolev inequalities we conclude
\[
|\p q |^2_{L^2}  \le |\CD q |^2_{L^2}  + (|q|^2_{L^4} | q |_{L^2})^2 \le C(E_0) |\CD q |_{L^2}^2
\]
where the Sobolev inequality was used in the last step with $\p q$ replaced by $\CD q$ which is true due to the observation
\[
|\p_k |q|^2| = 2 |\< q, \CD_k q\>| \le 2|q| |\CD_k q| \in L^1.
\]

To obtain $H^1$ bounds on $q$ multiply equation \eqref{hypc}
by $\overline{\CD_t q_\ell }$ and take the real part to obtain
\begin{equation}\label{h2}
\frac 12 \frac d{dt} \int  |\CD_k q_\ell |^2
+ \int \<q_k, q_0\> \< \CD_k q_\ell, i  q _\ell \>
+ \<  q_\ell, iq_k \> \< \CD_t q _\ell, iq_k \>
+ \varepsilon  |\CD_t q_\ell |^2 = 0 .
\end{equation}
Writing $\CD$ for the spatial covariant derivative, the second term can be bounded by
$$\begin{aligned}
\int |  \<q_k, q_0\> \< \CD_k q_\ell, i  q _\ell \>| dx \le & | \CD_k q_\ell|_{L^2}|q_0|_{L^6}|q|^2_{L^6}
\le | \CD_k q_\ell|_{L^2}|q_0|_{L^2}^{\frac 13}|\CD q_0|_{L^2}^{\frac 23}|q|_{L^2}^{\frac 23}|\CD q|_{L^2}^{\frac 43}\\*
 \le&C(E_0)  | \CD q|_{L^2}^{\frac 73}|q_0|_{L^2}^{\frac 13}|\CD q_0|_{L^2}^{\frac 23}\\*
 \le&C(E_0)  |\CD q|_{L^2}^{\frac 83}|\CD q_0|_{L^2}^{\frac 23} \qquad {\rm since} \ (\varepsilon - i )q_0 = \CD_j q_j\\*
 \le&\frac1{\sqrt{\varepsilon}}C(E_0)  |\CD q|_{L^2}^{4} + \frac {\varepsilon}4|\CD q_0|_{L^2}^2,
\end{aligned}
$$
and the third term by
$$
\int| \<  q_\ell, iq_k \> \< \CD_t q _\ell, iq_k \>|dx \le C |q|^3_{L^6}|\CD_t q _\ell|_{L^2}\le \frac{C(E_0)}{\varepsilon} |\CD q|_{L^2}^{4} + \frac {\varepsilon}4|\CD q_0|_{L^2}^2
$$
Using the identities $\CD_k q_0 = \CD_t q_k$, the above inequality, and equation \eqref{h2}, we have
\[
\frac 12 \frac d{dt} \int  |\CD_k q_\ell |^2
+ \frac \varepsilon 2 \int | \CD_t q _\ell|^2
\le \frac{C(E_0)}{\varepsilon} |\CD q|_{L^2}^{2} \int  |\CD_k q_\ell |^2 .
\]
Since by the energy identity  \( \varepsilon \int^\infty_0 |\CD q(t)|^2_{L^2} \,dt \le E_0 \)
we obtain global bounds on $\CD q$
\[
\int  |\CD_k q_\ell (t)  |^2
 + \varepsilon \int^t_0 \int | \CD_t q _\ell|^2
 \le C e^{\frac {C(E_0)}{\varepsilon^2}},
\]
which implies the desired bound on $\p q$.

Using this smooth solution $q$ we can construct a global smooth solution $u$ by means of proposition \ref{p2}. To show that $u$ belongs to the stated spaces we only need  to show that
\[
\int (u_0 -1) dx \le C.
\]
In fact, equation \eqref{pert} is equivalent to
\[
\p_t u = \frac 1{1+ \varepsilon^2} ( \eta \partial_k (u \times \partial_k u ) +\varepsilon (\Delta u - \< \nabla u , \eta \nabla u \> u)).
\]
Integrating the first component we obtain
\[
\frac d {dt} \int (u_0 -1 )
 = - \frac \varepsilon{1+ \varepsilon^2}\int \| \nabla u \|^2  \cosh \chi dx \le 0.
\]
\end{proof}
Weak solutions to \sms\ into $\h^2$ can be constructed as weak limits of the above solutions as $\varepsilon\to0$.
\begin{theorem}\label{page12}
Given \( u_* \in \dot H ^1 (\RR^2, \h^2 ) \)
such that
\( \int (u_{*_0} -1 ) dx < \infty \)
there exists a global weak solution to the \sms\ system
\[
\frac{\partial u}{\partial t} = \eta \partial_\ell (u \times \partial_\ell u )
 \qquad  u(0) = u_*
\]
with \( \|\nabla u\| \in L^\infty ( L^2 (\RR^2) ),\
u \in C (L^2_{\rm loc}  (\RR^2) )\
\text{and}\
u_0 -1 \in C (L^1 (\RR^2 ) ) .\)
\end{theorem}

\begin{proof}
Approximate the initial data by smooth functions $u_{*k}$ so that $\| \nabla u_{*k} \|_{L^2 (\RR^2)}^2$ and $|u_{*k0} -1|_{L^1(\RR^2)}$ are uniformly bounded and the geodesic distance between $u_{*k}(x)$ and $u_*(x)$ on $\HH^2$ converges to $0$ in $L^2(\RR^2)$. Even though $\HH^2$ is not compact, this can still be done since $\HH^2$ is diffeomorphic to $\RR^2$. In fact, using hyperbolic coordinates $u_0= \cosh \chi$,
$u_1 = \sinh \chi \cos \theta$ and $u_2 = \sinh \chi \sin \theta$, one can first approximate $u_*$ by a map whose image is in a compact set and then modify it into a smooth map by standard methods. In the hyperbolic coordinates, the boundedness of $\| \nabla u_{*k} \|_{L^2 (\RR^2)}^2$ and $|u_{*k0} -1|_{L^1(\RR^2)}$ takes the form
\[
 \int ( |\nabla \chi_{*k}|^2 + \sinh^2 \chi_{*k} |\nabla \theta_{*k} |^2 ) dx
< C
\quad \text{and} \quad
\int ( \cosh \chi_{*k} -1 ) dx < C.
\]
From Proposition~\ref{dgl} we have a  global smooth
solution to
\begin{equation}\label{hypk}
\begin{split}
\partial_t u_k + \frac 1k \eta (u_k \times u_{k_t})
 &= \eta \partial_\ell (u_k \times \partial_\ell u_k )\\
u_k (0, x) & = u_{*k} (0, x )
\end{split}
\end{equation}
such that
\[ \begin{split}
\int \| \nabla u_k (t) \|^2 \, dx
& = \int | \nabla \chi_k (t) |^2
 + \sinh \chi_k (t) |\nabla \theta_k (t) |^2 \,dx
\le E_0\\
\int (u _{0k} (t) -1 ) dx
 & = \int ( \cosh \chi_k -1 ) dx
\le C.
\end{split}\]
Thus $\chi_k$ is bounded in $L^\infty (H^1)$ and by
Moser-Trudinger inequality $\forall$ compact $\Omega \subset \RR^2$
\[
\intl_\Omega \exp\left(\alpha \frac{\chi^2_k (t)}{E_0}\right) \, dx
\le C( \Omega)
\]
for some positive $\alpha$.  This implies as in Proposition \ref{lss} that for
a subsequence
\[ \begin{split}
\chi_k  &\rightharpoonup \chi \quad \text{weak} * \ \text{in}\ L^\infty(H^1)\\
u_k & \rightharpoonup u \quad  \text{weak} * \ \text{in}\ L^\infty(W^{1, p}_{{\rm loc}}), \qquad p\in [1, 2)
\end{split}\]
where  $\|\nabla u(t)\|_{L^2(\RR^2)} \le C$ and $\int (u_0(t)-1)dx < C$.
Moreover for every cut off function $\varphi \in C^\infty_0 (\RR^2 \times \RR)$
\[
\p_t (\varphi u_k)- \partial_t \varphi u_k
 + \frac 1k \eta (\varphi u_k \times u_{k_t})
 = \eta \partial_\ell (\varphi u_k \times \partial_\ell u_k )
- \eta \partial_\ell \varphi u_k \times \partial_\ell u_k
\]
which implies that $\varphi u_k$ is bounded in $H^1_{{\rm loc}} (\R, W^{-1, p} (\RR^2))$, $1\le p<2$.
Consequently we have a subsequence where
\[
u_k \to u \quad \text{in}\ C(L^p) \quad \text{locally and a.e.}
\]
These bounds allow us to pass to the limit in
equation \eqref{hypk} to obtain
\[ \begin{split}
&\partial_t u  = \eta \partial_\ell ( u \times \partial_\ell u )\\
& u (0, x)  = u_* (x)
\end{split} \]
in the sense of distributions.
\end{proof}

{\section{Epilogue} \label{Epilogue}}
The results stated in this paper can be generalized to compact Hermitian symmetric K\"ahler manifolds $(M,g,J)$.  The equivalence of the \sms\ system and the frame system can be done in an identical manner provided there exist global smooth Coulomb frames when the dimension of $M$ is greater than $2$. To show global existence of weak
solutions in any space dimension we need to write the \sm\  system in divergence form. Therefore we have to restrict ourselves to the case when  $M$ has vanishing  first cohomology group.
In such a setting one uses the Killing vector fields to define weak solutions to the {\sm} system (SM)
$$ \p_tu =JD_ku\p_ku.
$$
in the following manner:

 A vector field $X\in TM$ is called Killing if $\CL_X g =0$ and $\CL_X J=0$.  Consequently if one considers the one form $\omega$ defined by $\omega(V)= g(JX,V)$ then $\omega $ is a closed one form since $M$ is K\"ahler. Moreover  since the first cohomology vanishes $\omega$ is exact.  Whence there is a function $f_X$ such that $\omega = df_X$, and for a solution $u$ to the (SM) system we have
 $$
\begin{aligned}
 \p_t f_X(u) & =  \omega(u_t) = g(JX(u),u_t) = -\<X(u),D_k\p_k u\> \\
& = - \p_k \<X(u),\p_ku\> +  \<D_kX(u), \p_ku\> =  - \p_k \<X(u),\p_ku\>
\end{aligned}
$$
since $X$ is Killing.  If the $2n$-dimensional manifold $M$ is compact and has $m$ Killing vector fields $\{X_a\}_{a=1}^m$ such that $TM = \text{span}\{ X_1, \ldots, X_m\}$, then the (SM) system is equivalent to
$$\p_t f_{X_a }(u) =  - \p_k \<X_a (u),\p_ku\>, \qquad  a = 1,\cdots, m.
$$

\subsection*{Remarks} Though $\HH^2$ is not compact, actually the definition \eqref{E:smh} of weak solutions of \sms  \; targeted on $\HH^2$ can also be viewed in this formulation with two Killing vector fields $X_1 = J \nabla (\sinh \chi \cos \theta)$ and $X_2 = J\nabla (\sinh \chi \sin \theta)$.

Weak solutions in higher dimensions can also be constructed using the idea in \cite{sh88, sh97, fr96}.  In this case we 1) embed $M$ isometrically and equivariantly in $\RR^L$ \cite{ms80}, and 2) define $d(u)$ the  distance  function from $M$ to $u$ and let $\sigma >0$ be so that $d(u)$ is smooth in the tubular neighborhood $O = \{ u\in \RR^L \mid d(u) < \sigma \} $ of $M$. Extend $d$ globally as a smooth function
$$F(u) = \varphi(d) d +  (1- \varphi (d)) \sigma, \qquad d=d(u)
$$
where $\varphi \in C_0^\infty(-\sigma, \sigma)$ and $ 0\le \varphi \le 1$ and $\varphi|_{[-\frac \sigma2, \frac \sigma2]} =1$.  3) Extend $J$ smoothly to act on $T\RR^L$.  This can be achieved by first extending $J(p)$ for $p\in M$ to act on  $T\RR^L= T_pM \oplus T_pM^\perp$ by first projecting on $TM$ and then applying $J$.  This operator can be extended to $O$ as a constant in the directions normal to $M$, i.e., $\forall u\in O$  decompose $u =p + n$ where $p\in M$ and $n \perp T_pM$ and define $\tilde{J}(u) =J(p)$ acting on $T_uO$. Finally define $\hat{J}(u) = \varphi(d(u)) \tilde{J}(u)$ for $u\in \RR^L$. It is clear that $\hat J$ is skew-symmetric. 4) Solve the equation
\[
\begin{aligned}
\varepsilon\p_t^2 u  - \hat{J}(u)\p_t u   -\Delta u + \frac 1\delta F(u)F'(u) =0\\
u(0,x)= u_\ast(x)\in M, \qquad \p_t u(0,x) =0.
\end{aligned}
\]
which has conserved energy
$$
\int \varepsilon|\p_t u|^2 + | \nabla u|^2 +  \frac 1\delta |F(u)|^2 dx = \int | \nabla u_\ast|^2
$$
By the energy method, the above equation has global solutions in $H^1$. For any Killing vector field $X \in TM$, from the equivariance of the embedding, $X$ can be extended to a vector field $X: \RR^L \to T\RR^L$ which generates an isometry on $\RR^L$ and satisfies $X \perp F(u)F'(u)$. Therefore, we have
$$
\varepsilon  \p_t \<X(u),\p_tu\>   -\<X(u), \hat J(u)\p_t u\>   - \p_k \<X(u),\p_ku\> =0
$$
By letting $\delta\to 0$ we have from the energy identity $u\to M$ in the $L^2$ sense and whence the limit satisfies
$$
\varepsilon  \p_t \<X(u),\p_tu\>   + \p_t f_X(u)   - \p_k \<X(u),\p_ku\> =0.
$$
Finally as $\varepsilon\to 0$ we obtain the \sm\ system in conservation form.
\begin{nth} Given $u_\ast: \RR^d\to M$ such that $\nabla u_\ast \in L^2$, the \sm\ system
$$\begin{aligned}
 &\p_tu =JD_k\p_ku\\
 &u(0,x) = u_\ast(x),
 \end{aligned}
$$
has a global weak solution such that $u\in C(\RR,L^2)\cap L^\infty({\dot H}^1)$.
\end{nth}

\end{document}